\documentclass[onefignum,onetabnum]{siamart190516}

\usepackage{graphicx,float}
\usepackage{subfigure}
\usepackage{amsfonts}
\newsiamremark{remark}{Remark}
\newtheorem{assumption}[theorem]{Assumption}
\newcommand{\BF}{\mathcal{F}}

\newcommand{\R}{\mathbb{R}}

\newcommand{\tr}{\mbox{tr}}

\newcommand{\esssup}{\operatorname*{ess\,sup\:}}

\usepackage{hyperref}
\hypersetup{
colorlinks=true, %colorise les liens
breaklinks=true, %permet le retour ??? la ligne dans les liens trop longs
urlcolor= blue, %couleur des hyperliens
linkcolor= red, %couleur des liens internes
citecolor=blue, %couleur des r???f???rences
pdftitle={}, %informations apparaissant dans
pdfauthor={}, %les informations du document
pdfsubject={} %sous Acrobat.
}

\emergencystretch=\maxdimen
\title{State-Dependent Temperature Control for Langevin Diffusions\thanks{We thank the two reviewers and the Associate Editor for constructive comments which have led to an improved version of the paper.} }

\author{Xue feng
			Gao\thanks{Department of Systems
			Engineering and Engineering Management, The Chinese University of Hong Kong, Hong Kong, China.
			(\email{xfgao@se.cuhk.edu.hk}).}
	\and Zuo Quan Xu\thanks{Department of Applied Mathematics, The Hong Kong Polytechnic University, Hong Kong, China. (\email{maxu@polyu.edu.hk}). }
	 \and	Xun Yu Zhou\thanks{Department of Industrial Engineering and Operations Research and The Data Science Institute, Columbia University, New York, NY 10027, USA. (\email{xz2574@columbia.edu}). }
		}

\begin{document}

\maketitle

%% --------------------------------------------------------------------------------------------------------------
\begin{abstract}
We study the temperature control problem for Langevin
diffusions in the context of non-convex optimization. The classical
optimal control of such a problem is of the bang-bang type, which is
overly sensitive to errors. A remedy is to allow the diffusions to explore other temperature values and hence smooth out the bang-bang control. We accomplish this by a stochastic relaxed control formulation incorporating
randomization of the temperature control and regularizing its entropy. We derive a state-dependent, truncated exponential distribution, which
can be used to sample temperatures in a Langevin algorithm, in terms of the solution to an HJB partial differential equation.
We carry out a numerical experiment on a one-dimensional baseline example, in which the HJB equation can be easily solved, to compare the performance of the algorithm with three other available algorithms in search of a global optimum.

%\textbf{Kewords.}
%Langevin
%diffusion, non-convex optimization, stochastic relaxed control, entropy regularization, Boltzmann exploration, HJB equation.
%\end{keywords}

\end{abstract}


% REQUIRED

\section{Introduction}

Consider the problem of finding the global minimizer of a \textit{non-convex} function $f: \mathbb{R}^d \rightarrow \mathbb{R}$, where $f$ is assumed to be differentiable. Traditional algorithms such as gradient descent often end up at a local optimum. The simulated annealing (SA) technique \cite{kirkpatrick1983optimization} has been developed to resolve the problem
of algorithms being trapped in local optima. The main thrust of an SA algorithm is {\it exploration via randomization}: At each iteration,
the algorithm {\it randomly}
samples a solution close to the current one and moves to it according to
a (possibly time-dependent) probability distribution, which in the literature is mostly exogenous as a pre-defined schedule. This scheme allows for a more extensive search or exploration for the global optima with the risk of moving to worse solutions at some iterations. The risk is however controlled by
slowly cooling down over time the ``temperature" which is used to characterize the level of
exploration or randomization.

The {\it Langevin algorithm} applies the SA technique by adding an independent series of Gaussian noises to the classical gradient descent algorithm, where the variance
of the noises is linearly scaled by a sequence of temperature parameters $\{\beta_k\}$ which control the level of exploration/randomization. The continuous-time version of
the Langevin algorithm is the so-called \textit{overdamped} (or {\it first-order}) {\it Langevin diffusion} governed by a stochastic differential equation (SDE), where the temperature is now a stochastic process $\{\beta(t): t\geq0\}$. There is a large volume of literature on Langevin diffusions for studying non-convex optimization; see, to name just a few, \cite{chiang1987diffusion, geman1986diffusions, marquez1997convergence}. The advantage of studying a diffusion process instead of a discrete-time
iteration lies in the simplicity and tractability of the former thanks to the many available analytical tools such as stochastic calculus, stochastic control and partial differential equation (PDE).

If we fix the temperature process in a Langevin diffusion to be a constant $\beta>0$,
then one can prove that, under some mild assumptions on $f$,
the Langevin diffusion process converges to a unique stationary distribution whose density is a {\it Gibbs measure}, $\frac{1}{Z(\beta)}e^{- \frac{f(x)}{\beta}}$,
where $Z(\beta)$ is the normalizing factor \cite{chiang1987diffusion}. %Hence, Langevin diffusions can be used for sampling possibly high dimensional noises in an SA algorithm.
As the temperature cools down (i.e. as $\beta \rightarrow 0$), this stationary distribution increasingly concentrates around the global minimizers of $f$ which in turn justifies the use of the Langevin diffusion for non-convex optimization. Indeed, Langevin diffusions and their variants have recently found many applications in data science including large-scale Bayesian inference and stochastic non-convex optimization arising from machine learning; see e.g. \cite{dalalyan2017further, gurbuzbalaban2020decentralized, Raginsky17, welling2011bayesian} and the references therein.

%When $u_t = u$ is a positive constant for all $t$, then the diffusion process $X(\cdot)$ in \eqref{eq:overdamped} admits a unique stationary distribution with the density $\pi(x) \propto e^{-\frac{1}{u} f(x)}$, under some mild assumptions on $f$ (see e.g. \cite{chiang1987diffusion,stroock-langevin-spectrum}). For $u$ chosen properly (small enough), it is easy to see that this distribution will concentrate around approximate global minimizers of $f$. This provides an intuitive explanation for the success of Langevin diffusion for global optimization. Indeed, the Brownian noise helps the trajectory of the diffusion escape suboptimal local minimum of the non-convex objective $f$; Without the noise term, $X(\cdot)$ in \eqref{eq:overdamped} reduces to the gradient flow dynamics or gradient descent in continuous-time, and it is well-known that running gradient descent on a non-convex function may converge to a sub-optimal local minimum.

A critical question for solving global non-convex optimization using the Langevin diffusions is the design of the temperature schedule $\{\beta(t): t\geq0\}$.
In the literature, the temperature $\beta(\cdot)$ is typically taken either as a constant (a hyperparameter) or various functions of time $t$, mostly {\it exogenously} given; see e.g. \cite{chiang1987diffusion, gelfand1991recursive, geman1986diffusions, holley1989asymptotics}. When $\beta(t)\equiv \beta>0$, it is well known that the expected transition time between different local minima for the overdamped Langevin diffusion is exponential in $1/\beta$, a phenomenon known as metastability \cite{Bovier2004, bovier2005metastability}.
%\cite{Raginsky17} provide a non-asymptotic analysis of the time-discretized Langevin diffusion with stochastic gradients
%finite-time guarantees for SGLD to find approximate minimizers of both
%empirical and population risks.
%that the expected transition time
%is also exponential in the problem dimension $d$, based on which they carefully bound the discretization error when the objective is non-convex.
\cite{Raginsky17} and \cite{xu2018global} further upper bound the time of Langevin dynamics converging to an approximate global minimum of non-convex functions, and
\cite{zhang2017hitting} and \cite{chen2020stationary} analyze the hitting time of Langevin dynamics to a local minimum, where the temperature is all assumed to be a constant. On the other hand, \cite{munakata2001temperature} formulate a deterministic optimal control problem to investigate optimal temperature
schedules and derive an ordinary differential equation (ODE) for the time-dependent temperature $\beta(\cdot)$.

Ideally, in implementing an algorithm for finding the global optimizers of a function, the temperature should be fine tuned based on where the {\it current} iterate is; namely it should be a function of the state. For instance, in order to quickly escape a deeper local minimum, one should impose
a higher temperature on the process. On the other hand, only lower temperatures are needed when the current process is not trapped in local minima in order to stay near a good solution. As a result of the need for this state-dependence, the temperature should be formulated as a {\it stochastic process} because the state itself follows a stochastic process.
\cite{fang1997improved} consider a temperature process as a specific, \textit{exogenously} given increasing function of the current function value $f(x)$, so that the trajectory can
perform larger steps in the search space when the current solution is seen to be far from optimal.
They show numerically that this scheme yields more rapid convergence to global optimizers. %Other exogenously devised state-dependent temperatures have also been shown to improve the rate of convergence of the Langevin samplers to the target Gibbs measure using Markov chain Monte Carlo (MCMC); see

The goal of this paper is to develop a {\it theoretical} framework (instead of a heuristic approach) for designing
\textit{endogenous} state-dependent temperature schedule for Langevin diffusions in the context of non-convex optimization.
%One main contribution of this paper is to propose a novel formulation for the optimal temperature control problem, using (relaxed) stochastic control with entropy regularization.
It is natural to formulate an optimal stochastic control problem in which the controlled dynamic is the Langevin diffusion with the temperature process taken as
the control, while the objective is to minimize a cost functional related to the original function $f$ to be minimized.
However, because the temperature appears linearly in the variance term of the diffusion, optimal control is of a ``bang-bang" nature, i.e. the temperature ought to repeatedly switch between two extreme points which are the lower and upper bounds of the temperature parameter. Such a bang-bang solution is not desirable in practice because it is {\it extremely} sensitive to errors: a tiny error may cause
the control to take the {\it wrong} extreme point; see e.g. \cite{bertrand2002new, silva2010smooth}. To address this issue, we take the so-called {\it exploratory
control} formulation, developed by \cite{wang2019reinforcement}, in which
classical controls are randomized and replaced by {\it distributions} of controls.
This allows the control to take other values than just the extreme ones.
Moreover, to encourage certain level of randomization, we explicitly incorporate the entropy of the distributions - which measures the extent of randomization - into the cost functional as a {\it regularization} term.
This formulation smoothes out the temperature control and motivates the solution to deviate, if discreetly, from the overly rigid bang-bang strategy.

Entropy regularization is a widely used heuristic in reinforcement learning (RL) in discrete time and space; see
e.g. \cite{Haarnoja2018}. \cite{wang2019reinforcement} extend it to continuous time/space using a stochastic relaxed control formulation. The motivation behind \cite{wang2019reinforcement}'s formulation is {\it repeated learning} in unknown
environments. Policy evaluation is achieved by repeatedly sampling controls from a same distribution and applying law of large numbers. In the present paper there is not really an issue of {\it learning} (as in reinforcement {\it learning}) because we can assume that $f$ is a known function. The commonality that prompts us to use the same formulation of \cite{wang2019reinforcement} is {\it exploration} aiming at escaping from possible ``traps". Exploration is to broaden search in order to get rid of over-fitted solutions in RL, to jump out of local optima in SA, and to
deprive of too rigid schedules in temperature control.

Under the infinite time horizon setting, we show that
the optimal exploratory control is a truncated exponential distribution whose domain is the pre-specified range of the temperature. The parameter of this distribution is state dependent, which is determined by a nonlinear elliptic PDE in the general multi-dimensional case and an ODE in the one-dimensional case. The distribution is a continuous version of the {\it Boltzmann} (or {\it softmax}) {\it exploration},
a widely used heuristic strategy in the discrete-time RL literature \cite{Bridle90, Sutton18} that applies
exponential weighting of the state-action value function (a.k.a. {\it Q-function}) to balance exploration and exploitation. In our setting, the Q-function needs to be replaced by the (generalized) {\it Hamiltonian}. This, however, is not surprising because in classical stochastic control theory optimal feedback control is to {\it deterministically} maximize the Hamiltonian \cite{yong1999stochastic}.

%use the idea of entropy regularization \cite{Haarnoja2018, wang2019reinforcement} to randomize the control, where classical controls are replaced by distribution of controls. This smoothes out the temperature control and allows the solution to deviate from the bang-bang strategy.
%Mathematically, our formulation explicitly incoporates randomization of the temperature into the optimization objective as a regularization term, with a weight imposed on the differential entropy of the probability distributions of the temperature.
%We discuss this in details in Section~\ref{sec:control-entropy} and review the connection of the entropy-regularized optimization with overdamped Langevin diffusions in Section~\ref{sec:entropy-reg}.

%Boltzmann softmax operator can trade-off well between exploration and exploitation according to current estimation in an exponential weighting scheme
%hich can be solved numerically.

As discretization of the optimal state process in our exploratory framework, which satisfies an SDE, naturally leads to a Langevin algorithm with state-dependent noise,
our results have algorithmic implications whenever the HJB equations involved can be
easily solved. In particular, the HJB equation in our setting is an ODE when the state space is one-dimensional, which can be efficiently solved. For numerical demonstration,
we compare the performance of this new algorithm with three other existing algorithms for the global minimization of a one-dimensional baseline non-convex function. The first method is a Langevin algorithm with a constant temperature \cite{Raginsky17}, the second one is a Langevin algorithm where the temperature decays with time in a prescribed power-law form \cite{ neelakantan2015adding, welling2011bayesian}, and the last one is a replica exchange algorithm \cite{dong2020replica}.
 The experiment indicates that, at least in the one-dimensional case, our Langevin algorithm
with state-dependent noise finds the global minimizer faster and outperforms the other three methods.

We must, however, emphasize that the main contribution of this paper is {\it theoretical}, rather than {\it algorithmic}. It
establishes and develops a theoretical framework for studying state-dependent temperature control in SA, in which a non-convex optimization problem
is connected to an HJB PDE. Among other theoretical interests in this connection,
one observation is that solving the HJB requires {\it global} information on the underlying function $f$, whereas most existing SA
algorithms use only {\it local} information of $f$ at any given iteration. Intuitively, the former is more advantageous other things being equal. This is actually demonstrated by the outperformance of our algorithm in the one-dimensional case when solving the HJB equation poses no numerical challenge. The insight is that
one should try to make use of the global information of $f$ as much as possible.
That said, for high dimensional non-convex optimization, we are not advocating our theory for actual algorithmic implementation before the curse of dimensionality for PDEs has been resolved. Recently there has been some encouraging progress using deep neural networks to numerically and efficiently solve high-dimensional PDEs \cite{Beck2019,Han18}; so hopefully our results could also contribute to devising SA algorithms in the future when that line of research has come to full fruition.
%, so we expect that, in the future, solving an HJB equation may be an easier task than finding the global minimizer for a non-convex optimization problem. There are many other technics (such as fast Fourier transformation) on solving HJB equations may be applied simultaneously.
%Our second remark is that our algorithm requires global information to search for the global minimizer for non-convex optimization problem, by contrast, many algorithms using only local information to do that. On one hand it is a disadvantage as global information may not be available or it takes more time to implement. On the other hand, it is just unavoidable to require the global information in order to find the global minimizer; local information algorithms may waste time on searching for local minimizers as they do not use the global information properly. Fully taking the advantage of the global information may allow us to find the global minimizer easier and quicker than other methods.
%}

The rest of the paper proceeds as follows. In Section~\ref{sec:problem}, we describe the problem motivation and formulation. Section~\ref{sec:control} presents the optimal temperature control. In Section~\ref{sec:numerical}, we report numerical results comparing the performance of our algorithm with three other methods. Finally, we conclude in Section~\ref{sec:conclusion}.

%%%%%%%%%%%%%%%%%%%%%%%%%%%%%%%%%%%%%%%%%%%%%%%%%%%%%%%%%%%%%%%%%%%%%%
%%%%%%%%%%%%%%%%%%%%%%%%%%%%%%%%%%%%%%%%%%%%%%%%%%%%%%%%%%%%%%%%%%%%

\section{Problem Background and Formulation}\label{sec:problem}
%\textbf{Background.}
\subsection{Non-convex optimization and Langevin algorithm}
Consider a non-convex optimization problem:
\begin{equation}\label{original}
\min_{x\in \mathbb{R}^d} f(x) ,
\end{equation}
where $f$: $\mathbb{R}^{d}\rightarrow [0, \infty)$ is a continuously differentiable,
{\it non-convex} function. The Langevin algorithm aims to have global convergence guarantees and has the following iterative scheme:
\begin{equation} \label{eq:GLD-iterates}
X_{k+1}=X_{k}-\eta f_x(X_k) +\sqrt{2\eta\beta_k}\xi_k\,,\;\;k=0,1,2,\cdots,
\end{equation}
where $f_x$ is the gradient of $f$, $\eta>0$ is the step size, $\{\xi_k\}$ is i.i.d Gaussian noise and $\{\beta_k\}$ is a sequence of the temperature parameters (also referred to as a cooling or annealing schedule) that typically decays over time to zero.
This algorithm is based on the discretization of the
overdamped Langevin diffusion:
\begin{equation}\label{eq:overdamped}
dX(t)=-f_x(X(t))dt+\sqrt{2 \beta(t)}dW(t),\;\;X(0)=x \in \mathbb{R}^d,
%\qquad
%X(0)=x\in\mathbb{R}^{d},
\end{equation}
where $x$ is an initialization, $\{W(t):t\geq0\}$ is a standard $d$-dimensional Brownian motion with $W(0)=0$, and $\{\beta(t):t\geq0\}$ is an adapted, nonnegative stochastic process, both defined on a filtered probability space $(\Omega,{\cal F},\mathbb{P};\{\BF_t\}_{t\geq0})$ satisfying the usual conditions.

When $\beta(t)\equiv \beta>0$,
under some mild assumptions on $f$, the solution of \eqref{eq:overdamped} admits a unique stationary distribution with the density $\pi(x) \propto e^{-\frac{1}{\beta} f(x)}$. Moreover,
Raginsky et al. (2017) show that for a finite $\beta>0$,
\begin{equation}
	\mathbb{E}_{X \sim \pi }f (X) - \min_{x\in \mathbb{R}^d} f(x) \leq \mathcal{I}(\beta) := \frac{d \beta}{2 }\log\left(\frac{eM}{m}\left(\frac{b}{d \beta}+1\right)\right), \label{cor:const:I}
\end{equation}
where $\mathbb{E}_{X \sim \pi }f (X):= \int_{ \mathbb{R}^d} f(x) \pi(dx)$,
$M, m, b$ are constants associated with $f$. It is clear that $\mathcal{I}(\beta) \rightarrow 0$ when $\beta \rightarrow 0$. %, which indicates that the stationary distribution concentrates around the global minimizers.

Our problem is to
%
%\subsection{An Optimal Control Problem.}
%We take the continuous-time perspective and consider the following problem: given $T$ steps (the computational budget) to run \eqref{eq:overdamped-2}, how to
control the temperature process $\{\beta(t): t\ge 0\}$ so that the performance of the continuous Langevin algorithm \eqref{eq:overdamped} is optimized. We measure the performance using the expected total discounted values of the iterate $\{X(t): t\geq0\}$, which is $\mathbb{E}\int_0^\infty e^{-\rho t}f(X(t))dt$ where $\rho>0$ is a discount factor. The discounting puts smaller weights on function values that happen in later iterations; hence in effect it dictates the computational budget (i.e. the number of iterations budgeted) to run the algorithm \eqref{eq:overdamped}.
Clearly, if this performance or cost functional value (which is always nonnegative) is small, then it implies that, on average, the algorithm strives to achieve smaller values of $f$ over iterations, and terminates at a budgeted number of iterations, $T$, with an iterate $X(T)$ that is close to the global optimum.

Mathematically, given $\rho>0$, an arbitrary initialization $X(0)=x \in \mathbb{R}^d$ and the range of the temperature $U=[a,1]$ where {$0<a<1$}, we aim to solve the following stochastic control problem where
the temperature process is taken as the control:
\begin{eqnarray}\label{eq:opt1}
\begin{array}{ll}
\quad \quad  \mbox{Minimize} & \mathbb{E}\int_0^\infty e^{-\rho t}f(X(t))dt,\\
\quad  \quad  \mbox{subject to} & \left\{\begin{array}{l}
\mbox{equation (\ref{eq:overdamped})},\\
\beta:=\{\beta(t):t\geq 0\} \mbox{ is adapted, and } \beta(t)\in U \mbox{ a.e.} \;\
t\ge 0, \mbox{ a.s.}
\end{array}\right.
\end{array}
%_{(\beta_t)_{0 \le t \le T}: 1 \ge \beta_t \ge 0} \left[ \mathbb{E}_{x_0}f(X(T)) - \min_{x \in \mathbb{R}^d} f(x) \right],
\end{eqnarray}

%where $\mathbb{E}_{x_0}$ denotes the expectation conditioned on $X(0)=x_0$.
%Since $\min_{x \in \mathbb{R}^d} f(x)$ is a constant independent of the control $(\beta_t)_{0 \le t \le T}$, it is equivalent to solve
%\begin{equation} \label{eq:opt1}
%\min_{(\beta_t)_{0 \le t \le T}: 1 \ge \beta_t \ge 0} \mathbb{E}_{x_0}f(X(T)) .
%\end{equation}
%Here the constraint $\beta_t \ge 0$ can be replaced by $\beta_t >0$ if needed. It is expected that the optimal $\beta_t$ should depend on the location $X(t)$ and the geometry of $f$; for instance, in order to quickly escape a deeper local minima, the process $X(t)$ should have a high temperature / larger noise.
%An immediate follow-up question is how large $T$ should be (how many iterations are needed) so that $\mathbb{E}f(X_T) - \min_{x \in \mathbb{R}^d} f(x) \le \epsilon$ for a given accuracy level $\epsilon>0.$ It would be also interesting to study the same problem (and how to choose the stepsizes) for the discrete-time algorithm \eqref{eq:GLD-iterates} with stochastic gradients.

\begin{remark}[Choice of the range of the temperature]
{\rm Naturally the temperature has a nonnegative lower bound. We suppose that it also has an upper bound, assumed to be 1 without loss of generality. In the Langevin algorithm and SA literature, one usually uses a determinist temperature schedule $\beta(t)$ that is either a constant or decays with time $t$. Hence there is a natural upper bound which is the initial temperature. This quantity is generally problem dependent and can be a hyperparameter chosen by the user. The upper bound of the temperature should be tuned to be reasonably large for the Langevin algorithm to overcome all the potential barriers and to avoid early trapping into a bad local minimum.
}
%(e.g. on how deep the valleys/local minima are) in the study of theoretical convergence of SA algorithms; see e.g. \cite{Hajek88}. In addition, there have been quite a few heuristic/computational methods proposed to set the initial temperature of SA algorithms; see e.g. \cite{Ameur04} and the references therein for details.
%In modern machine learning applications, the initial temperature is usually a hyperparameter to be chosen by the user which is again problem dependent. For instance, a commonly used scheme in machine learning applications is $\beta(t) = \frac{c}{(a+t)^b}$ where $a, b, c>0$ are hyperparameters; see e.g. Ye et al. (2017). The upper bound on $\beta(\cdot)$ is an important parameter for our problem when it comes to numerical inplementation. It should be large enough for the diffusion to overcome all the potential/energy barriers and to avoid early trapping into a bad local minimum.
%Using a physical interpretation, it means the initial mean kinetic energy (or the temperature) of the particle/ball be high enough to overcome any energy barrier, therefore the ball can possibly pass through all the minima on its trajectory.
\end{remark}

%What if replace $\beta_t \in \mathbb{R}$ by a positive semidefinite matrix $\Sigma_t \in \mathbb{R}^{d \times d}$?

\subsection{Solving problem \eqref{eq:opt1} classically} \label{sec:bangbang}
Define the optimal value function of the problem (\ref{eq:opt1}):
\begin{eqnarray} \label{eq:opt-control-value}
\quad \quad  V_0(x)= \inf_{\beta \in \mathcal{A}_0(x)} \mathbb{E} \left[ \int_0^{\infty} e^{-\rho t} f(X(t)) dt \Big| X(0)=x \right],
%\text{subject to} \quad dX(t)=-\nabla f(X(t))dt+\sqrt{2 u_t}dW_{t}.
\end{eqnarray}
where $x\in \mathbb{R}^d,$ and 
$\mathcal{A}_0(x)$ is the set of admissible controls $\beta$ satisfying the constraint in (\ref{eq:opt1}). 
A standard dynamic programming argument \cite{yong1999stochastic} yields that $V_0$ 
satisfies the following Hamilton--Jacobi--Bellman (HJB) equation: %($u_t$ represents the control $\beta_t$)
\begin{equation} \label{eq:bang-PDE}
-\rho v(x) + f(x) + \min_{\beta \in [a,1]} \left[ \beta \tr(v_{xx}(x)) - f_x(x) \cdot v_x( x) \right] = 0, \quad x\in \mathbb{R}^d,
\end{equation}
where $\tr(A)$ denotes the trace of a square matrix $A$, and ``$x\cdot y$" the inner product between two vectors $x$ and $y$.
%with terminal condition $v(T, x ) = f(x).$

%\xu{Intuitively speaking, when $\tr(v_{xx}(x))<0$, the value function is of concave shaped near $x$, so its minimum shall be far away from $x$ and thus one should heat at the highest possible temperature $\beta^*(x)= 1$ to escape from this point $x$. Similarly when $\tr(v_{xx}(x))>0$, one should heat at the lowest possible temperature $\beta^*(x)=a$ to search for the possible minimum nearby. This is clearly reflected by the HJB equation. Indeed, }
%Then,
The standard verification theorem in stochastic control theory \cite[Chapter 5]{yong1999stochastic} yields that an optimal {\it feedback} control policy should achieve the minimum in the above equation. However, the term inside the min operator is linear in the control variable $\beta$; hence
the optimal policy
has the following bang-bang form:
$\beta^*(x)= 1$ if $\tr(v_{xx}(x))<0$, and $\beta^*(x)= a$ otherwise.
In economics terms, intuitively, the value $v(x)$ can be regarded as the disutility
of the resource $x$ (here, it is {\it dis}utility, instead of utility, because
the objective is to {\it minimize}). When $\tr(v_{xx}(x))<0$, $v$ is locally concave around $x$ suggesting a risk-seeking preference. Hence one should randomize at the maximum possible level.
A symmetric intuition applies to the opposite case when $\tr(v_{xx}(x))>0$.
This temperature control scheme shows that one should in some states heat at the highest possible temperature, while in other states cool down completely, depending on the sign of $\tr(v_{xx}(x))$. This is {\it theoretically} the optimal scheme; but practically
it is just too {\it rigid} to achieve good performance as it concentrates on two actions only,
thereby leaves no room for errors which are inevitable in computation. This motivates the introduction of relaxed controls and entropy regularization in order to {\it smooth out} the temperature process in the next subsection.

\begin{remark}[Bang-bang control and parallel/simulated tempering]
{\rm The bang-bang control policy bears some resemblance to the so-called parallel tempering (or replica exchange) and simulated tempering in the MCMC literature \cite{earl2005parallel, marinari1992simulated, tawn2020weight}, 
% \cite{earl2005parallel, swendsen1986replica, dupuis2012infinite} and simulated tempering \cite{marinari1992simulated, lee2018beyond}
although there are important differences. Both tempering algorithms swap processes with different temperatures based on certain mechanisms in order to sample the desired target distribution; so there is switching at work between
different pre-specified temperatures like the bang-bang control.
However, these algorithms are mostly designed for probabilistic {\it sampling} with Metropolis--Hastings rule used to accept or reject a swap between processes,
while our focus is (non-convex) {\it optimization}. Moreover,
the tempering algorithms use state space augmentation.
Specifically, parallel tempering
runs $n \ge 2$ copies of Langevin dynamics with different constant temperature parameters; hence the state space is $\mathbb{R}^{dn}$. Simulated tempering treats the temperature $\{\beta(t): t\ge 0\}$ as a stochastic process taking only finite values and augments it to the original dynamics so that the state space is $(d+1)-$dimensional. In contrast, if we apply the bang-bang control, we obtain only one copy of the Langevin dynamics where the temperature parameter is state-dependent which is not a constant.
Recently, \cite{dong2020replica} propose to use the replica exchange algorithm for non-convex optimization.
In the numerical experiments reported in Section~\ref{sec:numerical}, we will compare the performance of our algorithm with that of \cite{dong2020replica}. }
\end{remark}

% in \cite{dong2020replica} (see Algorithm 1 there) runs one copy of gradient descent $(X_k)$ (i.e. Langevin dynamics with temperature $\beta=0$) and one copy of Langevin dynamics $(Y_k)$ in \eqref{eq:GLD-iterates} with constant temperature $\gamma>0$.
%It swaps the positions of $X_k$ and $Y_k$ if $f(X_k) > f(Y_k)$, and outputs $X_N$ at the terminal time $N$.

%\begin{remark}[A trivial example]
%{\rm Consider a special case where $f(x)=x^2$. In this case, a solution to the HJB equation \eqref{eq:bang-PDE} is given analytically by
%\begin{equation}
%v(x) = \left\{\begin{array}{ll}
%x^{-\rho/2},&x>0\\
%0,& x\leq 0
%\end{equation}
%We always have $v_{xx}>0$ for $x>0$ and hence the optimal $\beta^*(x)=0$, $x>0$. That is, the Langevin algorithm becomes simply the gradient descent without noise, provided that one initialize at some $x_0>0$. This is not surprising since $f$ has only one global minimum without any local minima, and adding noise to the iterates is not necessary.}
%\end{remark}
%\textcolor{red}{Xuefeng: Maybe remove this trivial example. It requires $x>0$ which results in a constrained optimization problem. In \eqref{eq:bang-PDE} it is unconstrained. Langevin algorithm for constrained optimziation is not well studied, though there are some recent results}

\subsection{Relaxed control and entropy regularization}

We now present our entropy-regularized relaxed control formulation of the problem.
Instead of a classical control $\{ \beta(t): t \ge 0 \}$ where $\beta(t) \in U=[a,1]$ for $t\ge0$, we consider a {\it relaxed} control $ \pi=\{\pi(t,\cdot): t \ge 0\}$, which is a {\it distribution} (or randomization) of classical controls over the control space $U$ where a temperature $\beta(t) \in U$ can be sampled from this distribution whose probability density function is $\pi(t,\cdot)$ at time $t$.

Specifically, we study the following {\it entropy-regularized stochastic relaxed control problem}, also termed as {\it exploratory control problem}, following \cite{wang2019reinforcement}:
\begin{align} \label{eq:opt-control-new}
 V_\lambda(x)&= \inf_{\pi \in \mathcal{A}(x)} \mathbb{E} \Big[ \int_0^{\infty} e^{-\rho t} f(X^{\pi}(t)) dt  \\
& \quad \quad \quad \quad \quad \quad - \lambda \int_0^{\infty} e^{-\rho t} \int_{U} - \pi(t,u) \ln \pi(t,u) du dt \Big| X^{\pi}(0)=x \Big], \nonumber
%\text{subject to} \quad dX(t)=-\nabla f(X(t))dt+\sqrt{2 u_t}dW_{t}.
\end{align}
where $x\in \mathbb{R}^d,$ the term $\int_{U} - \pi(t,u) \ln \pi(t,u) du$ is the differential entropy of the probability density function $\pi(t,\cdot)$ of the randomized temperature at time $t$, $\lambda > 0$ is a weighting parameter for entropy regularization,
$\mathcal{A}(x)$ is the set of admissible distributional controls to be precisely defined below, and the state process $X^{\pi}$ is governed by
\begin{equation} \label{eq:inf-dynamics}
dX^{\pi}(t)=- f_x(X^{\pi}(t))dt+\tilde \sigma(\pi(t))dW(t),\;\;X^{\pi}(0)=x \in \mathbb{R}^d,
\end{equation}
where
\begin{equation} \label{eq:sigma}
\tilde{\sigma} (\mu): = \sqrt{\int_U 2u \mu(u) du},\;\;\mu\in \mathcal{P}(U),
\end{equation}
with $\mathcal{P}(U)$ being the set of probability density functions over $U$.

The form of $\tilde{\sigma}(\mu)$ follows from the general formulation presented in
\cite{wang2019reinforcement}.
The function $V_{\lambda}$ is called the (optimal) {\it value function}.

Denote by  $\mathfrak{B}(U)$ the Borel algebra on $U$.
\begin{definition}
{\rm We say a
density-function-valued stochastic process $\pi= \{\pi(t,\cdot): t \ge 0\}$, defined on a filtered probability space $(\Omega,{\cal F},\mathbb{P};\{\BF_t\}_{t\geq0})$ along with a standard $d$-dimensional $\{\BF_t\}_{t\geq0}$-adapted Brownian motion
$\{W(t):t\geq0\}$,
is an {\it admissible (distributional) control}, denoted simply by $\pi \in \mathcal{A}(x)$, if
\begin{itemize}
\item[(i)] For each $t\ge 0$, $\pi(t,\cdot) \in \mathcal{P}(U)$ a.s.;

\item[(ii)] For any $A\in \mathfrak{B}(U)$, the process $(t,\omega)\mapsto \int_A \pi(t,u,\omega)du$ is $\BF_t$-progressively measurable;
%\item $\sup_{(t,\omega)}\int_U 2|u| \pi_t(u,\omega) du<\infty$. This is automatically true since $U$ is a compact set.
\item[(iii)] The SDE \eqref{eq:inf-dynamics} admits solutions on the same filtered probability space, whose distributions are all identical.
\end{itemize}}
\end{definition}

So an admissible distributional control is defined in the {\it weak} sense, namely,
the filtered probability space and the Brownian motion are also {\it part} of the control. Condition (iii) ensures that the performance measure in our problem (\ref{eq:opt-control-new}) is well defined even under the weak formulation. For notational simplicity, however, we will henceforth only write $\pi \in \mathcal{A}(x)$.

\begin{remark}[Choice of the weighting parameter $\lambda$]
{\rm The weighting parameter $\lambda$ dictates a trade-off between optimization (exploitation) and randomization (exploration), and is also called a temperature parameter in \cite{wang2019reinforcement}. \cite{TZZ2021} prove through a PDE argument that $V_\lambda$ converges to $V_0$ as $\lambda\rightarrow0$ and provide a convergence rate. Note that this temperature $\lambda$ is different from the temperature $\beta$ in the original Langevin diffusion.
These are temperatures at {\it different} levels: $\beta$ is at a lower, local level related to injecting Gaussian noises to the gradient descent algorithm, whereas $\lambda$ is at a higher, global level associated with adding (not necessarily Gaussian) noises to the bang-bang control. It is certainly an open interesting question to endogenize $\lambda$, which is yet another temperature control problem at a different level. In implementing the Langevin algorithm, however,
$\lambda$ can be fine tuned as a hyperparameter. }
\end{remark}

\section{Optimal Temperature Control}\label{sec:control}

\subsection{A formal derivation}

To solve the control problem \eqref{eq:opt-control-new}, we apply Bellman's principle of optimality to derive the following HJB
equation:

\begin{eqnarray}\label{eq:ODE2}
\quad \quad \quad - \rho v(x) - f_x(x)\cdot v_x(x)+f(x)+\inf_{\pi\in\mathcal{P}(U) } \int_{U} \left[ \tr(v_{xx}(x))u+\lambda \ln \pi(u) \right] \pi(u) du =0.
\end{eqnarray}

The verification approach then yields that the optimal distributional control can be obtained by solving the minimization problem in $\pi$ in the HJB equation. This problem
has a constraint $\pi\in\mathcal{P}(U)$, which amounts to
\begin{equation}\label{nor}
\int_U \pi(u) du =1, \quad \pi(u) \ge 0 \quad \text{a.e. on $U$}.
\end{equation}
For any $k\in\R$,
\begin{align*}
\int_{U} \left[ \tr(v_{xx}(x))u+\lambda \ln \pi(u) \right] \pi(u) du
&=\int_{U} \left[ \tr(v_{xx}(x))u+\lambda \ln \pi(u) +k\right] \pi(u) du -k.
\end{align*}
The integrand on the right hand side is a convex function of $\pi(u)$; so the first-order condition gives its unique global minimum at
\begin{eqnarray} \label{eq:opt-control-01}
\bar{\pi}(u; x)
&=& \exp\left(-\frac{1}{\lambda}[\tr(v_{xx}(x))u ] -\frac{k}{\lambda} -1 \right) \nonumber\\
&:=&
\frac{1}{Z(\lambda,v_{xx}(x))}\exp\left(-\frac{1}{\lambda}[\tr(v_{xx}(x))u ] \right), \quad u \in U,\;\;x\in \mathbb{R}^d.
\end{eqnarray}
Since $\bar{\pi}(\cdot; x)$ must satisfy (\ref{nor}), we deduce
%where $Z(\lambda,v_{xx}(x))$ (which depends on $k$) is the normalizing factor that makes $\bar{\pi}(\cdot;x)$ a density function, with
\begin{equation*}
Z(\lambda,v_{xx}(x))=\int_U \exp\left(-\frac{1}{\lambda}[ \tr(v_{xx}(x))u ] \right) du>0.
\end{equation*}
%while $\mathbb{S}^d$ is the space of $d\times d$ symmetric matrices.
This gives the optimal {\it feedback law} $\bar{\pi}(\cdot;x)$, which
is a {\it truncated} (in $U$) {\it exponential distribution} with the (state-dependent) parameter $c(x):=\frac{\tr(v_{xx}(x))}{\lambda}$. Note we do not require $\tr(v_{xx}(x))>0$ (i.e. $v$ is in general non-convex) or $c(x)>0$ here.

Substituting \eqref{eq:opt-control-01} back to the \eqref{eq:ODE2}, and noting
\begin{eqnarray*}
\lefteqn{\int_{U} \left[ \tr(v_{xx}(x))u+\lambda \ln \bar{\pi}(u; x) \right] \bar{\pi} (u; x) du} \\
&= &\int_{U} \left[ \tr(v_{xx}(x))u - \tr(v_{xx}(x))u - \lambda \ln (Z(\lambda,v_{xx}(x))) \right] \bar{\pi} (u; x) du \\
&=& - \lambda \ln (Z(\lambda,v_{xx}(x))),
\end{eqnarray*}
we obtain
%\begin{equation}
%v_t(t, x) + \min_{\pi \in \mathcal{P}[0,1]} \int_{0}^1 \left[ u \cdot v_{xx}(t,x) - f'(x) v_x(t, x) + \lambda \ln \pi(u) \right] \pi(u) du = 0,
%\end{equation}
% leads to
the following equivalent form of the HJB equation, a nonlinear elliptic PDE:
% \begin{eqnarray*}
%&-\rho v(x) - f_x(x)\cdot v_x(x) +f(x)+ v_{xx}(t,x) \cdot \frac{\lambda}{ v_{xx}(t,x)} \left( \frac{1 - (\frac{v_{xx}(t, x)}{\lambda} +1) e^{-\frac{v_{xx}(t, x)}{\lambda} }}{ 1 -e^{- \frac{v_{xx}(t, x)}{\lambda} }} \right) \\
%&+ \lambda \left( \ln\left(\frac{\frac{v_{xx}(t, x)}{\lambda}}{1-e^{-\frac{v_{xx}(t, x)}{\lambda} }}\right) - \frac{1 - (\frac{v_{xx}(t, x)}{\lambda} +1) e^{-\frac{v_{xx}(t, x)}{\lambda} }}{ 1 -e^{-\frac{v_{xx}(t, x)}{\lambda} }} \right)=0,
%\end{eqnarray*}
%or equivalently,
\begin{eqnarray}\label{eq: PDE-constraint}
-\rho v(x) - f_x(x)\cdot v_x(x) +f(x) - \lambda \ln (Z(\lambda,v_{xx}(x))) =0, \quad x \in \mathbb{R}^d.
\end{eqnarray}

%Note that the term inside the logarithmic function in the above is well defined and positive regardless of the definiteness of the Hessian matrix $v_{xx}(x).$

\begin{remark}[Gibbs Measure and Boltzmann Exploration]
{\rm %The optimal feedback law in \eqref{eq:opt-control-01} is closely related to Boltzman (or softmax) exploration, which is a
In reinforcement leaning there is a widely used {\it heuristic} exploration strategy called the {\it Boltzmann exploration}.
In the setting of \textit{maximizing} cumulative rewards,
Boltzmann exploration uses the Boltzmann distribution (or the Gibbs measure) to assign
a probability $p(s_t, a)$ to action $a$ when in state $s_t$ at time $t$:
\begin{equation} \label{eq:boltzman}
p(s_t, a) = \frac{e^{Q_t(s_t, a)/ \lambda}}{ \sum_{a=1}^m e^{Q_t(s_t, a)/ \lambda }}, \quad a=1, 2, \ldots, m,
\end{equation}
where
$Q_t(s, a)$ is the Q-function value of a state-action pair $(s,a)$,
and $\lambda>0$ is a parameter that controls the level of exploration; see e.g. \cite{Bridle90, Cesa17, Sutton18}. When $\lambda$ is high, the actions are chosen in almost equal probabilities. When $\lambda$ is low, the highest-valued actions are more likely to be chosen, giving rise to something resembling the so-called {\it epsilon-greedy policies} in multi-armed bandit problems. When $\lambda =0$, (\ref{eq:boltzman}) degenerates into the Dirac measure concentrating on the action that
maximizes the Q-function; namely the best action is always chosen. In the setting of continuous time and continuous state/control, the Q-function $Q_t(s,a)$ is not well defined and can not be used to rank and select actions \cite{Tallec19}. However, we can use the (generalized) {\it Hamiltonian} \cite{yong1999stochastic}, which in the special case of the classical problem \eqref{eq:opt1} is
\begin{align*}
H(x, u, v_x, v_{xx}):= - f_x(x)\cdot v_x(x)+f(x)+ \tr(v_{xx}(x))u.
\end{align*}
Then it is easy to see that \eqref{eq:opt-control-01} is equivalent to the following:
\begin{eqnarray}\label{BE}
\bar{\pi}(u; x)
&=& \frac{\exp(-\frac{1}{\lambda}[ H(x, u, v_x, v_{xx}) ] )}{\int_{u \in U} \exp(-\frac{1}{\lambda}[ H(x, u, v_x, v_{xx}) ] ) du }, \quad u \in U,\quad x \in \mathbb{R}^d,
\end{eqnarray}
which is analogous to $\eqref{eq:boltzman}$ with $x$ being the state and $u$ being the action. The negative sign in the above expression is due to the minimization problem instead of a maximization one. The importance of this observation is that here we lay a {\it theoretical underpinning} of the Boltzmann exploration via
entropy regularization, thereby provide interpretability/explainability of a largely heuristic approach.\footnote{The formula (\ref{BE}) was derived in
\cite[eq. (17)]{wang2019reinforcement} for a more general setting, but the connection with Boltzmann exploration and Gibbs measure was not noted there.} }
\end{remark}

%The term softmax for the action selection rule \eqref{eq:boltzman} is due to \cite{Bridle90}. This rule appears to be based on Luce's choice axiom \cite{Luce59}. See \cite{Sutton18}.
%
%Boltzman (or softmax) exploration is a widely used heuristic exploration strategy, but as pointed out in \cite{Cesa17} there is virtually no theoretical understanding about the limitations or the actual benefits of this exploration scheme in reinforcement learning (RL).
%
%Several prior studies (see e.g. \cite{Odono16, Haarnoja17, Ziebart2010} and references therein) have explored the connection of Boltzman exploration with maximum entropy policies in discrete-time RL problems, where $Q$ values are obtained (approximated) from energy-based models such as deep neural networks and restricted Boltzman machines. Our problem is of continous-time, the definition of $Q$ is different (nontrivial), and it involves a PDE.

\subsection{Admissibility and optimality}

The derivation above is formal, and is legitimate only after we {\it rigorously} establish the verification theorem and check the admissibility of the distributional controls
induced by the feedback law (\ref{eq:opt-control-01}).

To this end we need to impose some assumptions.
%From PDE perspective, unlike the unconstrained case with $\pi_t(u)$ supported on $(0,\infty)$, the term inside $\ln$ now is well defined and positive regardless of the sign of $v_{xx}.$ In addition, it is not necessary to have $v_{xx} >0$ and the optimal control $\pi^*(u; t,x)$ in \eqref{eq:opt-control-01} can still make sense even when $v_{xx}<0$ due to the constraint on $u$ to the bounded interval $[0,1]$.

\begin{assumption}\label{boundedf'}
{\rm We assume
\begin{itemize}
\item[(i)] $U=[a,1]$ for some $0<a<1$.
\item[(ii)] $f\in C^1(\mathbb{R}^d)$ and $|f_x(x)|<C$ for some constant $C>0$.
%\item[(ii)] $f\in C^1$ and there exist constants $C>0$ and $\delta \in (0,1)$ such that $|f_x(x)|\le C (1+ |x|^{\delta})$ for all $x$.
\end{itemize}
}
\end{assumption}

Assumption (i) imposes some minimal heating, $a>0$, to the Langevin diffusion, capturing the fact that in practice we would generally not know whether we have reached the global minimum at any given time; therefore we always carry out exploration, however small it might be, until some prescribed stopping criterion is met. On the other hand, the difficulty of the original non-convex optimization lies in the possibility of many {\it local} optima, while it is relatively easier to know (sometime from the specific context of an applied problem under consideration) the approximate range of the possible location of a global optimum. Therefore, the assumption of the gradient of $f$ in Assumption (ii) is reasonable as the function value outside of the above range is irrelevant.

Apply the feedback law (\ref{eq:opt-control-01}) to the controlled system (\ref{eq:inf-dynamics}) to obtain the following SDE
\begin{equation}\label{feedbackSDE}
dX^{*}(t)=- f_x(X^{*}(t))dt+h(X^{*}(t))dW_{t},\;\;X^*(0)=x
%\qquad
%X(0)=x\in\mathbb{R}^{d},
\end{equation}
where
\begin{equation} \label{eq:h}
h(x):=\sqrt{\frac{2}{Z(\lambda,v_{xx}(x))}\int_a^1 u \exp\left(-\frac{1}{\lambda}[\tr(v_{xx}(x))u ] \right) du}.
\end{equation}

%\begin{proposition}\label{ws}
%The SDE (\ref{feedbackSDE}) has a unique weak solution.
%\end{proposition}
%
%\begin{proof}
It is clear that $0<\sqrt{2a}\leq h(x)\leq \sqrt{2}$. Thus the drift and diffusion coefficients of (\ref{feedbackSDE}) are both uniformly bounded, while the latter satisfies the uniform elliptic condition. It follows from \cite[p. 87, Theorem 1]{Kry80} that (\ref{feedbackSDE}) has a weak solution. Moreover, \cite[Theorem 6.2]{SV69} asserts that the weak solution is unique.
%\end{proof}

Now define for $t\ge0, u\in[a,1]$,
\begin{equation}\label{pistar}
\pi^*(t,u):= \bar{\pi}(u; X^{*}(t))
\equiv \frac{1}{Z(\lambda,v_{xx}(X^{*}(t)))}\exp\left(-\frac{1}{\lambda}[\tr(v_{xx}(X^{*}(t)))u ] \right),
\end{equation}
where $\{X^*(t):t\ge0\}$ is a weak solution to (\ref{feedbackSDE}). Let
$\pi^*:=\{\pi^*(t,\cdot):t\ge0\}$. Then clearly $\pi^*\in \mathcal{A}(x)$.

\begin{theorem}\label{verification}
Suppose $v\in C^2$ is a solution to the HJB equation \eqref{eq:ODE2}, and $|v(x)|<C(1+|x|^k)$ $\forall x \in \mathbb{R}^d$ for some constants $C>0$ and $k>0$. Then $v= V_\lambda$. Moreover, in this case
the distributional control $\pi^*$ defined by (\ref{pistar}) is optimal.
\end{theorem}

%\begin{remark}
%In the special case that $f(x)= |x|$ and $\lambda=1$, the solution $v(t,x ) \equiv |x|$ satisfies the PDE (except at location point 0). In this case, the optimal distributional control $\pi^*(u; t,x)$ in \eqref{eq:opt-control-01} becomes a uniform distribution on $[0,1]$, and it is independent of $(t,x ).$ \gao{If there are more special cases for which the solution $v$ and the optimal control $\pi^*$ can be explicitly solved, that would be useful.}
%\end{remark}

We need a series of technical preliminaries to prove this theorem.

\begin{lemma}\label{entropyestimate}
For any $ \pi \in \mathcal{A}(x)$, we have
\begin{align}
\int_{a}^1 \pi(t,u) \ln \pi(t,u) du \geq a, \;\mbox{a.s.}
\end{align}
%Here $|U|$ is the Lebesgue measure of $U$.
\end{lemma}
\noindent
{\bf Proof}.
Applying the general inequality $x\ln x\ge x-1$ for $x> 0$, we have
\begin{align}
\int_{a}^1 \pi(t,u) \ln \pi(t,u) du&\geq \int_{a}^1 (\pi(t,u)-1)du= a.
\end{align}
%\end{proof}

\begin{lemma}\label{momentsestimates}
Suppose $ \{X(t): t \ge 0\}$ follows
\[dX(t)=b(t)dt+\sigma(t)dW({t}),\;\;t\ge0\]
with
\[C_b:=\esssup_{(t, \omega)} |b(t)|<\infty,\quad C_{\sigma}:=\esssup_{(t, \omega)} |\sigma(t)|<\infty. \]
Then, for any $k\ge 1$, there exists a constant $C>0$, which is independent of $T$ and $X(0)$, such that
\begin{equation}\label{1st}
\mathbb{E} \left[\sup_{0\le t\leq T}|X(t)|^{k}\right]<C(1+T^{k}+|X(0)|^{k}),\;\;\forall T\ge 0.\end{equation}
%Moreover,
%\begin{equation}\label{2nd}
%\mathbb{E} \left[e^{\varepsilon X(t)}\right]
%\leq e^{\varepsilon X(0)} \exp\bigg\{\Big(\varepsilon C_b+\frac{1}{2}\varepsilon^2 C_{\sigma}^2\Big)t\bigg\},\;\;\forall t\ge0, \;\varepsilon\ge 0.
%\end{equation}

\end{lemma}
\noindent
{\bf Proof}.
By the elementary inequality
\[(a+b+c)^{k}\leq (3\max\{a,b,c\})^k\leq 3^ka^k+3^kb^k+3^kc^k,\quad a,b,c\geq 0,\]
we have
\begin{align*}
\mathbb{E} \left[\sup_{0\le t\leq T} |X(t)|^{k}\right]&\leq \mathbb{E} \left[ \left(|X(0)|+\sup_{0\le t\leq T}\int_0^t|b(s)|ds+\sup_{0\le t\leq T}\left|\int_0^t\sigma(s)dW_{s}\right|\right)^{k}\right] \\
&\leq \mathbb{E} \left[ \left(|X(0)|+C_bT+\sup_{0\le t\leq T}\left|\int_0^t\sigma(s)dW_{s}\right|\right)^{k}\right] \\
&\leq 3^{k}|X(0)|^{k}+3^{k}C_b^{k}T^{k}+3^{k}\mathbb{E} \left[\sup_{0\le t\leq T}\bigg|\int_0^t\sigma(s)dW_{s}\bigg|^{k}\right]\\
&\leq 3^{k}|X(0)|^{k}+3^{k}C_b^{k}T^{k}+3^{k}C_k\left(\mathbb{E} \left[\int_0^T|\sigma(s)|^2ds\right]\right)^{k/2}\\
&\leq 3^{k}|X(0)|^{k}+3^{k}C_b^{k}T^{k}+3^{k}C_kC_{\sigma}^{k}T^{k/2},
\end{align*}
where the second to last inequality is due to the Burkholder-Davis-Gundy inequality.
This proves (\ref{1st}).

%Next, by It\^{o}'s lemma,
%\begin{align*}
%e^{\varepsilon X(t)}= e^{\varepsilon X(0)}+\int_0^t \Big[\varepsilon b(s)+\frac{1}{2}\varepsilon^2\sigma^2(s) \Big]e^{\varepsilon X(s)}ds+\varepsilon e^{\varepsilon X(s)}\sigma(s)dW_{s}.
%\end{align*}
%Let $\tau_n=\inf\{s\ge 0: |X(s)|>n\}$, $n\geq0$.
%It follows that
%\begin{align*}
%\mathbb{E} \left[e^{\varepsilon X(t\wedge \tau_n)}\right]
%&\leq e^{\varepsilon X(0)}+\mathbb{E} \left[\int_0^{t\wedge \tau_n}\Big[\varepsilon |b(s)|+\frac{1}{2}\varepsilon^2\sigma^2(s) \Big]e^{\varepsilon X(s)}ds \right]\\
%&\leq e^{\varepsilon X(0)}+\Big(\varepsilon C_b+\frac{1}{2}\varepsilon^2 C_{\sigma}^2\Big)\mathbb{E} \left[\int_0^{t\wedge \tau_n}e^{\varepsilon X(s)}ds \right]\\
%&\leq e^{\varepsilon X(0)}+\Big(\varepsilon C_b+\frac{1}{2}\varepsilon^2 C_{\sigma}^2\Big)\int_0^t\mathbb{E} \left[e^{\varepsilon X(s)} \right]ds.
%\end{align*}
%Sending $n\to \infty$, we have by Fatou's lemma,
%\begin{align*}
%\mathbb{E} \left[e^{\varepsilon X(t)}\right]
%&\leq e^{\varepsilon X(0)}+\Big(\varepsilon C_b+\frac{1}{2}\varepsilon^2 C_{\sigma}^2\Big)\int_0^t\mathbb{E} \left[e^{\varepsilon X(s)} \right]ds.
%\end{align*}
%Gronwall's inequality then yields (\ref{2nd}).
%\begin{align*}
%\mathbb{E} \left[e^{\varepsilon X(t)}\right]
%&\leq e^{\varepsilon X(0)} \exp\bigg\{\Big(\varepsilon C_b+\frac{1}{2}\varepsilon^2 C_{\sigma}^2\Big)t\bigg\}.
%\end{align*}
%\end{proof}

\begin{proposition}\label{V}
There exists a constant $C_1>0$ such that the value function $V$ satisfies
\begin{equation}
|V_\lambda(x)|\leq C_1(1+|x|), \;\;\forall x \in \mathbb{R}^d.
\end{equation}
\end{proposition}
\noindent
{\bf Proof}.
%\begin{proof}
It follows from Assumption \ref{boundedf'}-(ii) that $|f(x)|\leq C|x|+|f(0)|$.
Let $\pi$ be the uniform distribution on $[a,1]$.
It then follows from Lemma \ref{momentsestimates} that
\begin{align*}
V_\lambda(x) &\leq \mathbb{E} \left[ \int_0^{\infty} e^{-\rho t} (C|X^{\pi}(t)|+|f(0)|) dt+\lambda \int_0^{\infty} e^{-\rho t} \int_{a}^1 \pi(u) \ln \pi(u) du dt\right]\\
& <C'(|x|+1)
\end{align*}
for some constant $C'$ independent of $x$. On the other hand, for any $ \pi \in \mathcal{A}(x)$, by virtue of Lemma \ref{entropyestimate} we have
\begin{align*}
&\quad\;\mathbb{E} \left[ \int_0^{\infty} e^{-\rho t} f(X^{\pi}(t)) dt - \lambda \int_0^{\infty} e^{-\rho t} \int_{a}^1 - \pi_t(u) \ln \pi_t(u) du dt \right]\\
& \quad \quad \geq \mathbb{E} \left[ \int_0^{\infty} -e^{-\rho t}(C|X^{\pi}(t)|+|f(0)|) dt+\lambda \int_0^{\infty} e^{-\rho t} a dt\right] >-C''(|x|+1)
\end{align*}
for some constant $C''$ independent of $x$.
The proof is complete.

We have indeed established in the above that
\begin{align*}
\mathbb{E} \left[ \int_0^{\infty} e^{-\rho t} |f(X^{\pi}(t))| dt \right]<\infty.
\end{align*}
To solve problem \eqref{eq:opt-control-new}, we only need to consider those admissible controls $\pi$ such that
\begin{equation}\label{exclude}
\mathbb{E} \left[ \int_0^{\infty} e^{-\rho t} \int_{U} \pi(t,u) \ln \pi(t,u) du dt \right] <\infty,
\end{equation}
because any control violating (\ref{exclude}) renders the infinite value of the minimization problem \eqref{eq:opt-control-new} and hence can be excluded from consideration.
In view of this, (\ref{exclude}) is henceforth
assumed for any {\it admissible} distributional control.

%Next, we apply the feedback control (\ref{eq:opt-control-01}) to the controlled system (\ref{eq:overdamped-2}) to obtain the following SDE
%\begin{equation}\label{feedbackSDE}
%dX^{*}(t)=- f_x(X^{*}(t))dt+h(X^{*}(t))dW_{t},\;\;X(0)=x
%%\qquad
%%X(0)=x\in\mathbb{R}^{d},
%\end{equation}
%where
%\[ h(x):=\sqrt{\int_a^1 2u \frac{\exp(-\frac{1}{\lambda}[\tr(v_{xx}(x))u ] ) }{ \int_U \exp(-\frac{1}{\lambda}[ \tr(v_{xx}(x))u ] ) du}du}.
%\]
%
%\begin{proposition}\label{ws}
%The SDE (\ref{feedbackSDE}) has a unique weak solution.
%\end{proposition}
%
%\begin{proof}
%It is clear that $\sqrt{2a}\leq h(x)\leq \sqrt{2}$. Thus the drift and diffusion coefficients of \ref{feedbackSDE}) are both uniformly bounded, while the latter satisfies the uniform elliptic condition. It follows from \cite[p. 87, Theorem 1]{Kry80} that (\ref{feedbackSDE}) has a weak solution. Moreover, \cite[Theorem 6.2]{SV69} asserts that the weak solution is unique.
%\end{proof}
%
%Now define (with a slight abuse of notation)
%\[ \pi^*(t;u):= \frac{\exp(-\frac{1}{\lambda}[\tr(v_{xx}(X^*(t)))u ] ) }{ \int_U \exp(-\frac{1}{\lambda}[ \tr(v_{xx}(x))u ] ) du},\;\;t\ge0,\;u\in[a,1], \]
%where $\{X^*(t):t\ge0\}$ is a weak solution to (\ref{feedbackSDE}). Let
%$\pi^*:=\{\pi^*(t;\cdot):t\ge0\}$.
%
%\begin{proposition}
%\[\pi^*\in \mathcal{A}(x).\]
%\end{proposition}
%
%\begin{proof} This follows readily from Proposition \ref{ws}.
%\end{proof}

\bigskip

Now we prove Theorem \ref{verification}.

For any $ \pi \in \mathcal{A}(x)$, let $\{X^\pi(t):t\ge0\}$ be the solution to (\ref{eq:inf-dynamics}). Define $\tau_n:=\inf\{s\ge 0: |X^{\pi}(s)|>n\}$.
For any $T>0$, by It\^{o}'s lemma,
\begin{align*}
&\quad\;e^{-\rho(\tau_n\wedge T)}v(X^{\pi}(\tau_n\wedge T))\\
&=v(X^{\pi}(0))+\int_0^{\tau_n\wedge T} e^{-\rho t}v_x(X^{\pi}(t))\cdot \tilde{\sigma}(\pi(t))d W(t)\\
&+\int_0^{\tau_n\wedge T} e^{-\rho t}\Big(-\rho v(X^{\pi}(t))-f_x(X^{\pi}(t))\cdot v_x(X^{\pi}(t))+\frac{1}{2}{\tilde\sigma}^2(\pi(t))\tr(v_{xx}(X^{\pi}(t)))\Big)dt\\
&=\int_0^{\tau_n\wedge T} e^{-\rho t} v_x(X^{\pi}(t))\cdot \tilde{\sigma}(\pi(t))d W(t)
+\int_0^{\tau_n\wedge T} e^{-\rho t}\bigg( \int_a^1 \tr( v_{xx}(X^{\pi}(t))) u\pi(t,u) du
\\
&\quad\;-f(X^{\pi}(t))-\inf_{\pi'\in\mathcal{P}([a,1]) } \int_a^1 \left[ \tr( v_{xx}(X^{\pi}(t)))u+\lambda \ln \pi'(u) \right] \pi'(u) du\bigg)dt + v(x)\\
&\geq v(x)+\int_0^{\tau_n\wedge T} e^{-\rho t} v_x(X^{\pi}(t))\cdot \tilde{\sigma}(\pi(t))d W(t)\\
&\quad\;+\int_0^{\tau_n\wedge T} e^{-\rho t}\left( -f(X^{\pi}(t)-\lambda\int_a^1 \pi(t,u) \ln \pi(t,u) du )\right)dt.
\end{align*}
It follows
\begin{align*}
v(x)&\leq \mathbb{E} \left[e^{-\rho (\tau_n\wedge T) }v(X^{\pi}(\tau_n\wedge T))\right] \\
& \quad +\mathbb{E} \left[\int_0^{\tau_n\wedge T} e^{-\rho t}\left(f(X^{\pi}(t))+\lambda \int_a^1 \pi(t,u) \ln \pi(t,u) du\right)dt\right]\\
%&= \mathbb{E} \left[e^{-r (\tau_n\wedge T) }v(X^{\pi}(\tau_n\wedge T))\right]
%+\mathbb{E} \left[\int_0^{\tau_n\wedge T} e^{-\rho t}f(X^{\pi}(t))dt\right]\\
%&\quad\;+\mathbb{E} \left[\int_0^{\tau_n\wedge T} e^{-\rho t}\left(\lambda\int_a^1 \big( \pi(t,u) \ln \pi(t,u)-\pi(t,u)+1\big) du+\lambda a\right)dt\right]\\
&= \mathbb{E} \left[e^{-\rho (\tau_n\wedge T) }v(X^{\pi}(\tau_n\wedge T))\right]
+\mathbb{E} \left[\int_0^{\tau_n\wedge T} e^{-\rho t}f(X^{\pi}(t))dt\right]\\
&\quad\;+\lambda\mathbb{E} \left[\int_0^{\tau_n\wedge T} e^{-\rho t}\int_a^1 \pi(t,u) \ln \pi(t,u) dudt\right].%+\lambda a\mathbb{E} \left[\int_0^{\tau_n\wedge T} e^{-\rho t}dt\right].
\end{align*}
Sending $n\to \infty$, we have by Lemma \ref{momentsestimates} and the dominated convergence theorem that the sum of the first two terms converges to
\[ \mathbb{E} \left[e^{-\rho T}v(X^{\pi}(T))\right]+\mathbb{E} \left[\int_0^{T} e^{-\rho t}f(X^{\pi}(t))dt\right],\]
while it follows from (\ref{exclude})
along with the monotone convergence theorem (noticing $\pi(t,u) \ln \pi(t,u)\geq \pi(t,u)-1$) that the last term converges to
\begin{align*}
\lambda\mathbb{E} \left[\int_0^{T} e^{-\rho t}\int_a^1 \pi(t,u) \ln \pi(t,u) dudt\right].
\end{align*}
So
\begin{align*}
v(x)&\leq \mathbb{E} \left[e^{-\rho T}v(X^{\pi}( T))\right]
+\mathbb{E} \left[\int_0^{T} e^{-\rho t}f(X^{\pi}(t))dt\right]\\
&\quad\;+\lambda\mathbb{E} \left[\int_0^{T} e^{-\rho t}\int_a^1 \pi(t,u) \ln \pi(t,u) dudt\right].
\end{align*}
Noting Lemma \ref{momentsestimates} along with the assumption that $v$ is of polynomial growth, we have
\begin{equation}\label{bc}
\lim_{T\to\infty}\mathbb{E} \left[e^{-\rho T}v(X^{\pi}( T))\right]=0.\end{equation}
Applying the same limiting argument as above we let $T\rightarrow \infty$ to get
\[
v(x)\leq \mathbb{E} \left[\int_0^{\infty} e^{-\rho t}\left(f(X^{\pi}(t))+\lambda\int_U \pi(t,u) \ln \pi(t,u) du\right)dt\right].
\]
Since $\pi$ is arbitrarily chosen, it follows $v\leq V_\lambda$; in other words $v$ is a lower bound of $V_\lambda$.

All the inequalities in the above analysis become equalities when we take
$\pi=\pi^*\in \mathcal{A}(x)$ defined by (\ref{pistar}), because
$\pi^*$ attains the infimum in the HJB equation (\ref{eq:ODE2}). This yields the aforementioned lower bound is achieved by $\pi^*$; hence $\pi^*$ is optimal and $v=V_\lambda$.

\begin{remark}[Existence and Uniqueness of Solution to HJB]
{\rm Due to the linear growth of the value function $V_\lambda$ (Proposition \ref{V}), we can prove by a standard argument (part of it is similar to the argument used in the above proof) that $V_\lambda$ is a solution to the HJB equation \eqref{eq:ODE2} as long as $V_\lambda\in C^2$. Moreover, Theorem \ref{verification} yields that
the solution is unique among $C^2$ functions of polynomial growth. The uniqueness follows from (\ref{bc}) which is essentially a boundary condition of \eqref{eq:ODE2}. When $V_\lambda$ is not $C^2$, we can apply viscosity solution theory to study the well-posedness of the HJB equation. However, we will not pursue this direction as it is not directly related to the objective of this paper.}
\end{remark}

\section{Algorithm and Numerical Example} \label{sec:numerical}

As discussed in Introduction, the main contribution of this paper is to develop a theory that connects the state-dependent temperature control problem for SA to an HJB equation, the latter concerning the global information about the underlying function $f$. The theory also naturally gives rise to an algorithm that can be employed to solve the original non-convex optimization problem to the extent that the HJB equation can be efficiently solved. Specifically,
the optimal state process $\{X^{*}(t):t\ge 0\}$ is determined by the SDE \eqref{feedbackSDE}. Discretization of this equation then leads to a Langevin algorithm with state-dependent noise. One can use Euler-Maruyama discretization; see, e.g. \cite{mattingly2002ergodicity}, due to its simplicity and popularity in practice. We now report a numerical experiment
for a baseline non-convex optimization example in a one-dimensional state space for which the HJB equation is an ODE.

\begin{figure}[H]
\centering
\includegraphics[width= 0.7 \textwidth]{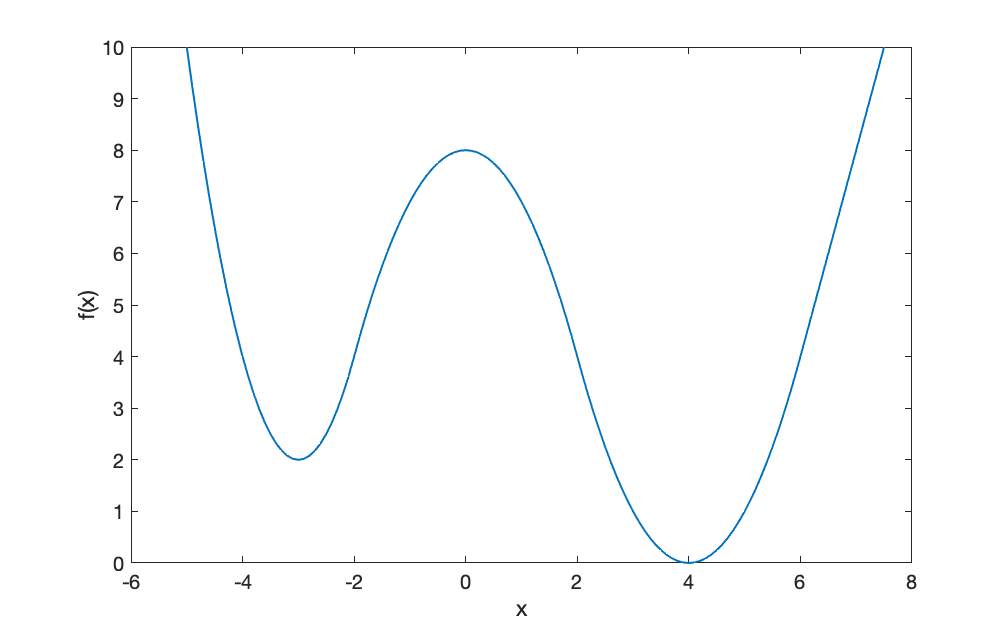}
\caption{A double well function $f$}
\label{fig:1}
\end{figure}

Consider the problem of finding the global minima of the following asymmetric double-well function:
\begin{equation}
f(x) = \begin{cases}
4x-20& \mbox{for} \quad x > 6, \\
	(x-4)^2 & \mbox{for} \quad 2< x \le 6 , \\
8 - x^2 & \mbox{for} \quad - 2 < x \leq 2 , \\
2(x+3)^2 +2 & \mbox{for} \quad -6 < x\leq -2 \\
-12 x - 52 & \mbox{for} \quad x\le -6.
\end{cases} \nonumber
\end{equation}
Note that $f: \mathbb{R} \rightarrow \mathbb{R}$ is a continuously differentiable function with two local minima at $ {x}_1 = -3$ and ${x}_2 = 4$, where the latter is the (unique) global minimum. See Figure~\ref{fig:1} for an illustration.

%\begin{figure}[h]
%\centering
%\includegraphics[width=0.5\textwidth]{doublewell}
% \caption{ODE solution $v(x)$ with initial condition given by $v(0)=-0.139, v'(0)=0.5026.$ }
%\label{fig:doublewell}
%\end{figure}

%It is clear that
%\begin{equation}
% f'(x) = \begin{cases}
%x - 1 & \mbox{for} \quad x \geq \frac{1}{2}, \\
% -x & \mbox{for} \quad -\frac{1}{8} \leq x \leq \frac{1}{2}, \\
% x+\frac{1}{4} & \mbox{for} \quad x\leq -\frac{1}{8}.
%\end{cases} \nonumber
%\end{equation}

We now present four algorithms, to be described below, for minimizing $f$ and compare their performances. All the algorithms are initialized at $x=-3$, the suboptimal local minimum of $f.$
We examine the convergence of the performance quantity $\mathbb{E}(f(X_k))$ where $X_k$ is the location of the iterate at iteration $k$, $k=1,2,\ldots$. In each of the algorithms, $\{\xi_k\}$ is a sequence of i.i.d standard normal random variable and $\eta$ is the constant step size which we tune respectively for best performance. In implementation, we change the upper bound of the allowed temperature to be 500, instead of 1 in our theoretical analysis, without loss of validity. This allows the iterates of the algorithms to climb over the barrier between the two local minimizers quickly.

The four algorithms are
\begin{itemize}
\item [(1)] Langevin algorithm with constant temperature $\beta$ \cite{Raginsky17}:
\begin{equation} \label{eq:GLD1-iterates}
X_{k+1}=X_{k}-\eta f'(X_k) +\sqrt{2\eta \beta} \xi_k\,.
\end{equation}
We choose the stepsize $\eta \in \{(1/2)^n: n =0,1,2, \ldots, 10 \}$ and the temperature $\beta \in \{ 500 \cdot (1/2)^n: n =0,1,2, \ldots, 15 \}$.
The best-tuned step size $\eta=(1/2)^1=0.5$ and temperature $\beta=500 \cdot (1/2)^{10}=0.4883.$

\item [(2)] Langevin algorithm with the power-law temperature schedule \cite{neelakantan2015adding}:
\begin{equation} \label{eq:GLD2-iterates}
X_{k+1}=X_{k}-\eta f'(X_k) +\sqrt{2\eta \beta_k} \xi_k\,,
\end{equation}
where we choose
$\beta_k = (\frac{d}{1+ k})^{b}$ with $b \in \{ 0.5, 0.6, \ldots, 1 \},$ and $d \in \{ 500\cdot (1/2)^n: n =0,1,2, \ldots, 7\}$; We also choose the stepsize $\eta \in \{(1/2)^n: n =0,1,2, \ldots, 10 \}$.
The best-tuned parameters are given by: $\eta =(1/2)^1=0.5,b =0.9, d=500\cdot (1/2)^4= 31.25$.

\item [(3)] Replica exchange (GDxLD in Algorithm 1 of \cite{dong2020replica}): one runs a copy of the gradient descent $\{X_k\}$ and a copy of the Langevin algorithm $\{Y_k\}$ with a constant temperature $\gamma>0$. If $f(X_k) > f(Y_k)$, then the positions of $X_k$ and $Y_k$ are swaped. Output $X_N$ as an optimizer of $f$ when the algorithm terminates at iteration $N.$ We choose $X_0=Y_0=-3$, and tune parameters where stepsize $\eta \in \{(1/2)^n: n =0,1,2, \ldots, 10 \}$, and constant temperature $\gamma \in \{ 500 \cdot (1/2)^n: n =0,1,2, \ldots, 15 \}$. The best-tuned parameters are given by $\eta =(1/2)^1=0.5, \gamma =500\cdot (1/2)^1= 250.$

\item [(4)] Our Langevin algorithm based on the Euler-Maruyama discretization of the SDE \eqref{feedbackSDE}:
\begin{equation} \label{eq:SDE-schedule}
X^{*}_{k+1}=X^{*}_{k}-\eta f'(X^{*}_k) +\sqrt{\eta} h( X^{*}_k) \xi_k\,.
\end{equation}
In dimension one we can deduce from \eqref{eq:h} that $h(x)=g \left(\frac{v_{xx}(x)}{\lambda}\right),$
where $v$ satisfies the ODE in \eqref{eq: PDE-constraint} with
\begin{eqnarray} \label{eq:g-oned}
Z(\lambda, v_{xx}(x))= \frac{ e^{- a \frac{v_{xx}(x)}{\lambda}} - e^{ -c \frac{v_{xx}(x)}{\lambda}}} { \frac{v_{xx}(x)}{\lambda}},
\end{eqnarray}
and
\begin{eqnarray} \label{eq:g-oned2}
g(y):=\sqrt{ 2 \frac{ (c+1/y) \cdot e^{-cy} - (a+1/y) \cdot e^{-y a}}{ e^{-cy} - e^{-y a}} }, \quad y \in (-\infty, \infty).
\end{eqnarray}
Here, $[a,c]$ is the range of the allowed temperature.
The hyperparameter $a$ in \eqref{eq:g-oned} is fixed at $0.0001$ for simplicity, and the hyperparameter $c$, which is the upper bound on the allowed temperature, is fixed at 500.
We choose $\rho, \lambda \in \{5\cdot (1/2)^n: n =0,1,2, \ldots, 8\}$, and stepsize $\eta \in \{(1/2)^n: n =0,1,2, \ldots, 10 \}$. To implement the algorithm \eqref{eq:SDE-schedule}, we need to solve the second-order nonlinear ODE in \eqref{eq: PDE-constraint}. We use random initial conditions $[v(0), v'(0)]$ where both $v(0)$ and $v'(0)$ are independently sampled from a standard normal distribution. We select 20 random initializations for the ODE.
The best-tuned parameters are given by $\eta=(1/2)^3=0.125, \lambda =5\cdot (1/2)^4=0.3125, \rho =5\cdot (1/2)^2=1.25$, with the ODE initialization $[v(0), v'(0)]=[-0.2853, 1.1575]$.

\end{itemize}

\begin{figure}[H]
			\centering
			\includegraphics[width= 0.75 \textwidth]{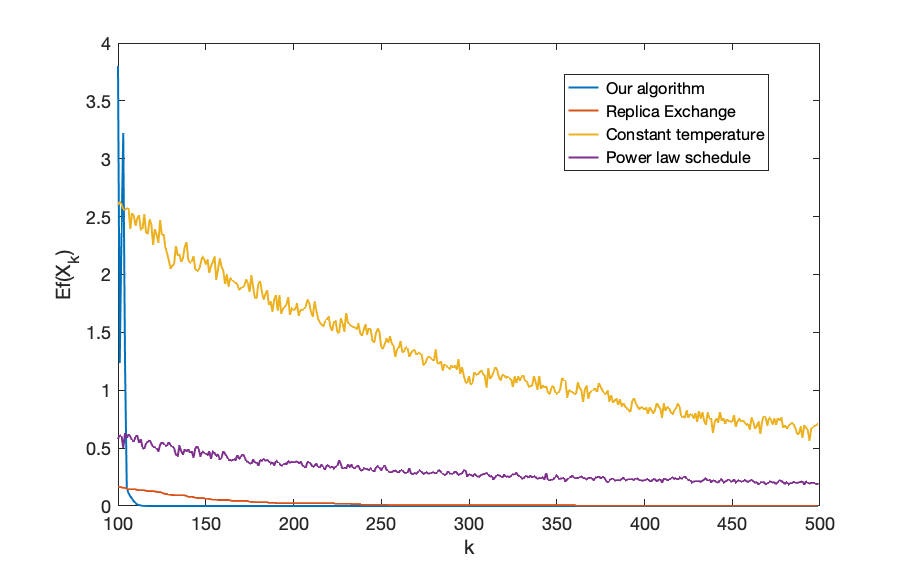}
			\caption{Performance comparison of the four algorithms, all initialized at $X_1=-3$.}
			\label{fig:comparison}
			\end{figure}

\begin{figure}[H]
			\centering
			\includegraphics[width= 0.75 \textwidth]{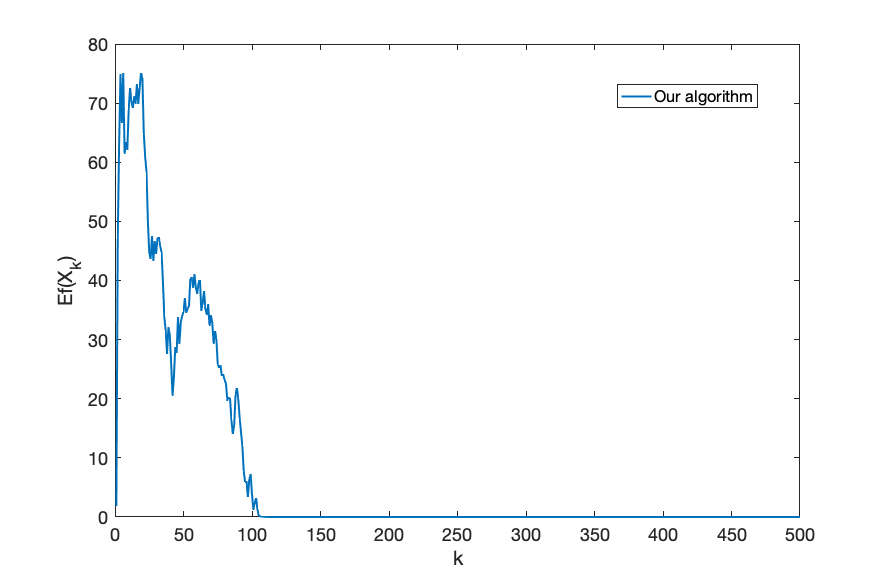}
			\caption{Performance of our algorithm initialized at $X_1=-3$.}
			\label{fig:ouralgo}
			\end{figure}

Figure \ref{fig:comparison} shows the performance of the four algorithms, where the expectation $\mathbb{E}(f(X_k))$
is approximated by its sample average with the sample size 500. Each algorithm is terminated if the number of iterations exceeds the allowed number of iterations, which we set to be 1000 in our experiment. In the first 100 iterations, the expected function values $\mathbb{E}(f(X_k))$ from our algorithm are very large, primarily because large noises are injected into the iterates during these initial iterations so as to escape the local minimum at $-3$; see Figure~\ref{fig:ouralgo} for a zoomed-in version. Hence,
for better visualization, we plot in Figure \ref{fig:comparison} the expected function values $\mathbb{E}(f(X_k))$ from $k=100$ up to only $k=500$ iterations.
As we can see, Langevin algorithms with constant temperature and power-law temperature schedule have difficulties in locating the global minima within 500 iterations. By its very definition, replica exchange is expected to perform better than the Langevin algorithm with a constant temperature. This is confirmed by Figure \ref{fig:comparison} in which
the replica exchange algorithm finds the global minimum quickly, although it needs to run two algorithms (the gradient descent and a Langevin algorithm) instead of only one algorithm.
Our Langevin algorithm with state-dependent noises can find the global minimum faster than the other three algorithms.
It is, however, computationally more expensive compared with other three methods,  needing to solve a nonlinear ODE. % at each iteration.

%\footnote{{\bf It will be better to have a more elaborate discussion here. Also, what happens to the sharp spike of the blue line near 0? Do the other two have similar spikes?}} %We also observe similar

% in around 25 steps, the algorithm \eqref{eq:SDE-schedule} converges to the global minimum.

\begin{figure}[H]
\centering
\subfigure[$v(x)$ on $(-6,4)$]{
\includegraphics[width=0.45\columnwidth]{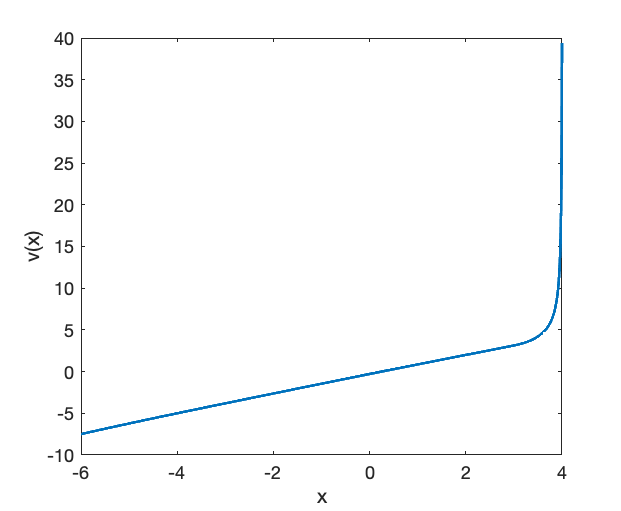}
\label{fig:vx}
}
\subfigure[$v''(x)$ on $(-6,2)$]{
\includegraphics[width=0.45\columnwidth]{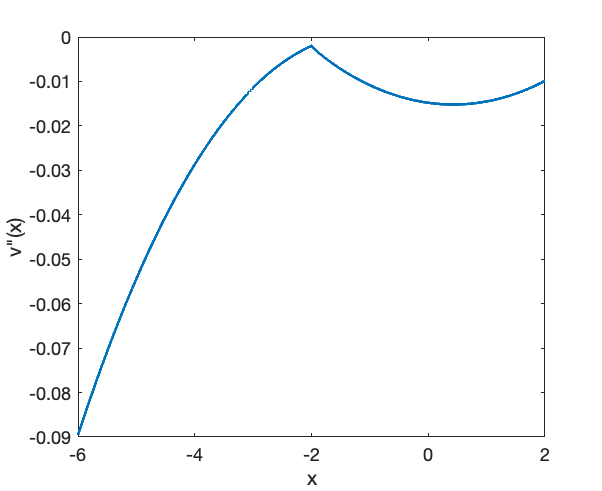}
\label{fig:vddx}
}
\caption{The HJB solution $v$ and its second order derivative $v''$.}
\end{figure}

\begin{figure}[H]
\centering
\includegraphics[width= 0.75\textwidth]{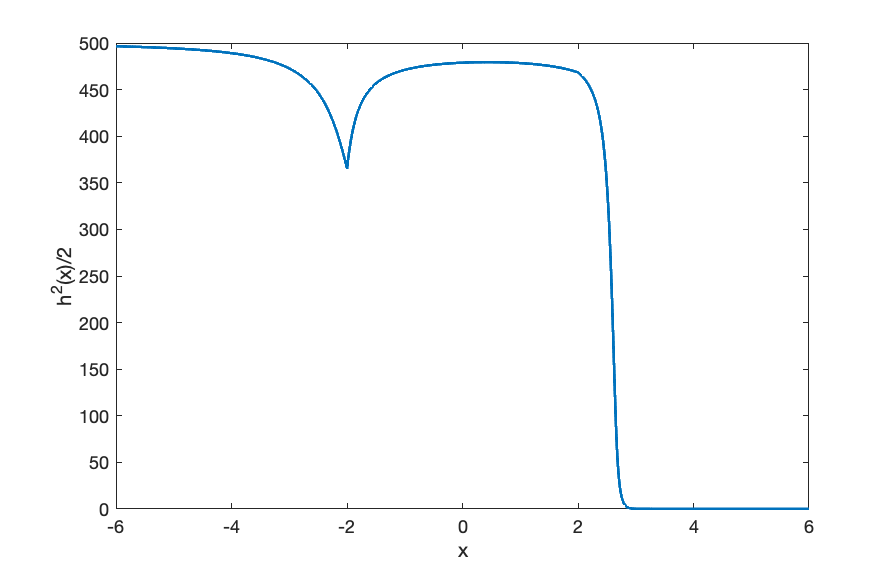}
\caption{State-dependent temperature function ${h^2}/2$.}
\label{fig:temperature2}
\end{figure}

The temperature process for our algorithm depends on the solution $v$ to the ODE \eqref{eq: PDE-constraint}. We plot in Figure~\ref{fig:vx} the function values $v(x)$ where $x\in (-6, 4)$. We observe that when $x$ is close to 4, the global minimum of $f$, $v$ grows very quickly. This is because much less exploration is needed near
the global minimum, in which case the subtraction of the entropy value from the overall objective value in (\ref{eq:opt-control-new}) is much less, boosting the value of $v$. We also plot
in Figure~\ref{fig:vddx} the second order derivative $v$ on the interval $(-6, 2)$.\footnote{We do not plot $v''(x)$ for $x>2$ in Figure~\ref{fig:vddx} for better visualization because $v''(x)$ becomes very large when $x>2$.}
However, $v$ or $v''$ affects the temperature process of our algorithm only indirectly, through the function $h$ given by \eqref{eq:SDE-schedule}.
So the plots of $v$ and $v''$ are less informative than that of the state-dependent temperature function $\frac{h^2}{2}$ which is proportional to the variance of the noise
in our algorithm; see Figure~\ref{fig:temperature2}. We see that the temperature is close to zero for $x>3$, and is mostly large elsewhere. This indicates that our state-dependent algorithm is ``intelligent": it uses the lowest temperature when close to the global minimum, and uses a large temperature (recall the largest temperature allowed is 500) to escape from traps such as suboptimal local minima or saddle points. We also observe that there is a prominent kink at $x=-2$ in Figure~\ref{fig:temperature2}. This is primarily due to the spike of $v''$ at $-2$ in Figure~\ref{fig:vddx} and the fact that $h(x)=g \left(\frac{v_{xx}(x)}{\lambda}\right)$ with $g$ given in \eqref{eq:g-oned2}.

\section{Conclusion}\label{sec:conclusion}

This paper aims to develop an endogenous temperature control scheme for applying Langevin
diffusions to find non-convex global minima. We take the exploratory stochastic control framework, originally proposed by \cite{wang2019reinforcement} for reinforcement learning, to account for the need of smoothing out the temperature process. We derive a state-dependent, Boltzmann-exploration type distributional control, which can be used to sample temperatures in a Langevin algorithm.
Numerical analysis shows that our algorithm outperforms three alternative ones based on
a constant temperature, a power decay schedule and a replica exchange method respectively. However, the function used in the numerical example is one-dimensional, for which the HJB equation is an ODE and hence easy to solve. For high-dimensional problems, the HJB equation is a PDE whose numerical solutions may suffer from the curse of dimensionality. Therefore, at least for now,
the main contribution of this paper is {\it not} algorithmic. Rather, it is, generally, to lay a {\it theoretical} underpinning for smoothing out often overly rigid classical optimal
controls (such as bang-bang controls) and, specifically, to provide an interpretable
state-dependent temperature control scheme for Langevin diffusions via HJB equations.\\

\section*{Acknowledgments}
Xuefeng Gao acknowledges financial support from Hong Kong GRF (No.14201117 and No.14201421).
Zuo Quan Xu acknowledges financial support from NSFC (No.11971409), Hong Kong GRF (No.15202817 and No.15202421), the PolyU-SDU Joint Research Center on Financial Mathematics and the CAS AMSS-POLYU Joint Laboratory of Applied Mathematics, The Hong Kong Polytechnic University.
Xun Yu Zhou acknowledges financial supports through a start-up grant at Columbia University and the Nie Center for Intelligent Asset Management. We also thank Mert G\"{u}rb\"{u}zbalaban and Lingjiong Zhu for their comments, and Yi Xiong for the help with the experiments.

\newpage

\bibliographystyle{siamplain}
\bibliography{SICON}

\begin{thebibliography}{10}

\bibitem{Beck2019}
{\sc C.~Beck, E.~Weinan, and A.~Jentzen}, {\em Machine learning approximation
  algorithms for high-dimensional fully nonlinear partial differential
  equations and second-order backward stochastic differential equations},
  Journal of Nonlinear Science, 29 (2019), pp.~1563--1619.

\bibitem{bertrand2002new}
{\sc R.~Bertrand and R.~Epenoy}, {\em New smoothing techniques for solving
  bang--bang optimal control problems -- numerical results and statistical
  interpretation}, Optimal Control Applications and Methods, 23 (2002),
  pp.~171--197.

\bibitem{Bovier2004}
{\sc A.~Bovier, V.~Gayrard, and M.~Klein}, {\em Metastability in reversible
  diffusion processes {I}: Sharp asymptotics for capacities and exit times},
  Journal of the European Mathematical Society, 6 (2004), pp.~399--424.

\bibitem{bovier2005metastability}
{\sc A.~Bovier, V.~Gayrard, and M.~Klein}, {\em Metastability in reversible
  diffusion processes {II}: Precise asymptotics for small eigenvalues}, Journal
  of the European Mathematical Society, 7 (2005), pp.~69--99.

\bibitem{Bridle90}
{\sc J.~S. Bridle}, {\em Training stochastic model recognition algorithms as
  networks can lead to maximum mutual information estimation of parameters}, in
  Advances in neural information processing systems, 1990, pp.~211--217.

\bibitem{Cesa17}
{\sc N.~Cesa-Bianchi, C.~Gentile, G.~Lugosi, and G.~Neu}, {\em Boltzmann
  exploration done right}, in Advances in neural information processing
  systems, 2017, pp.~6284--6293.

\bibitem{chen2020stationary}
{\sc X.~Chen, S.~S. Du, and X.~T. Tong}, {\em On stationary-point hitting time
  and ergodicity of stochastic gradient langevin dynamics.}, Journal of Machine
  Learning Research, 21 (2020), pp.~1--41.

\bibitem{chiang1987diffusion}
{\sc T.-S. Chiang, C.-R. Hwang, and S.~J. Sheu}, {\em Diffusion for global
  optimization in $\mathbb{R}^n$}, SIAM Journal on Control and Optimization, 25
  (1987), pp.~737--753.

\bibitem{dalalyan2017further}
{\sc A.~Dalalyan}, {\em Further and stronger analogy between sampling and
  optimization: Langevin monte carlo and gradient descent}, in Conference on
  Learning Theory, 2017, pp.~678--689.

\bibitem{dong2020replica}
{\sc J.~Dong and X.~T. Tong}, {\em Replica exchange for non-convex
  optimization}, arXiv preprint arXiv:2001.08356,  (2020).

\bibitem{earl2005parallel}
{\sc D.~J. Earl and M.~W. Deem}, {\em Parallel tempering: Theory, applications,
  and new perspectives}, Physical Chemistry Chemical Physics, 7 (2005),
  pp.~3910--3916.

\bibitem{fang1997improved}
{\sc H.~Fang, M.~Qian, and G.~Gong}, {\em An improved annealing method and its
  large-time behavior}, Stochastic processes and their applications, 71 (1997),
  pp.~55--74.

\bibitem{gelfand1991recursive}
{\sc S.~B. Gelfand and S.~K. Mitter}, {\em Recursive stochastic algorithms for
  global optimization in $\mathbb{R}^d$}, SIAM Journal on Control and
  Optimization, 29 (1991), pp.~999--1018.

\bibitem{geman1986diffusions}
{\sc S.~Geman and C.-R. Hwang}, {\em Diffusions for global optimization}, SIAM
  Journal on Control and Optimization, 24 (1986), pp.~1031--1043.

\bibitem{gurbuzbalaban2020decentralized}
{\sc M.~G{\"u}rb{\"u}zbalaban, X.~Gao, Y.~Hu, and L.~Zhu}, {\em Decentralized
  stochastic gradient langevin dynamics and hamiltonian monte carlo}, arXiv
  preprint arXiv:2007.00590,  (2020).

\bibitem{Haarnoja2018}
{\sc T.~Haarnoja, A.~Zhou, K.~Hartikainen, G.~Tucker, S.~Ha, J.~Tan, V.~Kumar,
  H.~Zhu, A.~Gupta, P.~Abbeel, et~al.}, {\em Soft actor-critic algorithms and
  applications}, arXiv preprint arXiv:1812.05905,  (2018).

\bibitem{Han18}
{\sc J.~Han, A.~Jentzen, and E.~Weinan}, {\em Solving high-dimensional partial
  differential equations using deep learning}, Proceedings of the National
  Academy of Sciences, 115 (2018), pp.~8505--8510.

\bibitem{holley1989asymptotics}
{\sc R.~A. Holley, S.~Kusuoka, and D.~W. Stroock}, {\em Asymptotics of the
  spectral gap with applications to the theory of simulated annealing}, Journal
  of functional analysis, 83 (1989), pp.~333--347.

\bibitem{kirkpatrick1983optimization}
{\sc S.~Kirkpatrick, C.~D. Gelatt, and M.~P. Vecchi}, {\em Optimization by
  simulated annealing}, Science, 220 (1983), pp.~671--680.

\bibitem{Kry80}
{\sc N.~Krylov}, {\em Controlled Diffusion Processes}, Springer, New York,
  1980.

\bibitem{marinari1992simulated}
{\sc E.~Marinari and G.~Parisi}, {\em Simulated tempering: a new monte carlo
  scheme}, EPL (Europhysics Letters), 19 (1992), p.~451.

\bibitem{marquez1997convergence}
{\sc D.~M{\'a}rquez}, {\em Convergence rates for annealing diffusion
  processes}, The Annals of Applied Probability,  (1997), pp.~1118--1139.

\bibitem{mattingly2002ergodicity}
{\sc J.~C. Mattingly, A.~M. Stuart, and D.~J. Higham}, {\em Ergodicity for
  {SDE}s and approximations: locally {L}ipschitz vector fields and degenerate
  noise}, Stochastic Processes and their Applications, 101 (2002),
  pp.~185--232.

\bibitem{munakata2001temperature}
{\sc T.~Munakata and Y.~Nakamura}, {\em Temperature control for simulated
  annealing}, Physical Review E, 64 (2001), p.~046127.

\bibitem{neelakantan2015adding}
{\sc A.~Neelakantan, L.~Vilnis, Q.~V. Le, I.~Sutskever, L.~Kaiser, K.~Kurach,
  and J.~Martens}, {\em Adding gradient noise improves learning for very deep
  networks}, arXiv preprint arXiv:1511.06807,  (2015).

\bibitem{Raginsky17}
{\sc M.~Raginsky, A.~Rakhlin, and M.~Telgarsky}, {\em Non-convex learning via
  stochastic gradient {L}angevin dynamics: A nonasymptotic analysis}, in
  Conference on Learning Theory, 2017, pp.~1674--1703.

\bibitem{silva2010smooth}
{\sc C.~Silva and E.~Tr{\'e}lat}, {\em Smooth regularization of bang-bang
  optimal control problems}, IEEE Transactions on Automatic Control, 55 (2010),
  pp.~2488--2499.

\bibitem{SV69}
{\sc D.~Stroock and S.~Varadhan}, {\em Diffusion processes with continuous
  coefficients, i}, Communications On Pure And Applied Mathematics, 22 (1969),
  pp.~345--400.

\bibitem{Sutton18}
{\sc R.~S. Sutton and A.~G. Barto}, {\em Reinforcement learning: An
  introduction}, MIT press, 2018.

\bibitem{Tallec19}
{\sc C.~Tallec, L.~Blier, and Y.~Ollivier}, {\em Making deep q-learning methods
  robust to time discretization}, arXiv preprint arXiv:1901.09732,  (2019).

\bibitem{TZZ2021}
{\sc W.~Tang, Y.~Zhang, and X.~Y. Zhou}, {\em The exploratory control problem
  with application to the state dependent temperature control for langevin
  diffusions}, working paper,  (2021).

\bibitem{tawn2020weight}
{\sc N.~G. Tawn, G.~O. Roberts, and J.~S. Rosenthal}, {\em Weight-preserving
  simulated tempering}, Statistics and Computing, 30 (2020), pp.~27--41.

\bibitem{wang2019reinforcement}
{\sc H.~Wang, T.~Zariphopoulou, and X.~Y. Zhou}, {\em Reinforcement learning in
  continuous time and space: A stochastic control approach}, Journal of Machine
  Learning Research,  (2019).

\bibitem{welling2011bayesian}
{\sc M.~Welling and Y.~W. Teh}, {\em Bayesian learning via stochastic gradient
  {L}angevin dynamics}, in Proceedings of the 28th International Conference on
  Machine Learning (ICML-11), 2011, pp.~681--688.

\bibitem{xu2018global}
{\sc P.~Xu, J.~Chen, D.~Zou, and Q.~Gu}, {\em Global convergence of {L}angevin
  dynamics based algorithms for nonconvex optimization}, in Advances in Neural
  Information Processing Systems, 2018, pp.~3122--3133.

\bibitem{yong1999stochastic}
{\sc J.~Yong and X.~Y. Zhou}, {\em Stochastic controls: Hamiltonian systems and
  HJB equations}, vol.~43, Springer Science \& Business Media, 1999.

\bibitem{zhang2017hitting}
{\sc Y.~Zhang, P.~Liang, and M.~Charikar}, {\em A hitting time analysis of
  stochastic gradient langevin dynamics}, in Conference on Learning Theory,
  2017, pp.~1980--2022.

\end{thebibliography}

\end{document}